\documentclass{amsart}

\usepackage[english]{babel}

\newtheorem{theorem}{Theorem}[section]
\newtheorem{definition}[theorem]{Definition}
\newtheorem{remark}{\it Remark\/}

\def\R{\mathbb R}
\def\N{\mathbb N}

\def\P{\mathbb P}
\def\ds{\displaystyle}

\begin{document}
\title{Finite time extinction for solutions to fast diffusion stochastic porous media equations}

\author{Viorel Barbu}
\email{vb41@uaic.ro}
\address{Institute of Mathematics ``Octav Mayer'', Iasi, Romania}

\author{Giuseppe Da Prato}
\email{daprato@sns.it}
\address{Scuola Normale Superiore di Pisa, Italy}

\author{Michael R\"ockner}
\email{roeckner@Mathematik.Uni-Bielefeld.DE}
\address{Faculty of Mathematics, University of Bielefeld, Germany and
Department of Mathematics and Statistics, Purdue University,\\  U. S. A}

\maketitle

\begin{abstract}
\selectlanguage{english}
{\it We prove that the solutions to
fast diffusion stochastic porous media equations have finite time extinction with strictly positive probability.}

\end{abstract}

\section{Introduction}
 Consider the stochastic porous media equation
  \begin{equation}
\label{e1.1}
\left\{\begin{array}{l}
dX(t)-\rho\Delta(|X|^\alpha(t)\;\mbox{\rm sign}\;X(t))dt-
\Delta(\tilde \Psi(X(t))dt=\sigma(X(t))dW(t),\;\mbox{\rm in}\; (0,\infty)\times  \mathcal O,\\
X=0\quad\mbox{\rm on}\; (0,\infty)\times \partial \mathcal O,\quad X(0,x)=x\quad\mbox{\rm on}\; \mathcal O, \\
\end{array}\right.
\end{equation}
where $\rho>0$, $\alpha\in (0,1)$, $\tilde \Psi$ is a   continuous monotonically non decreasing  function of linear growth and $\sigma(X)dW=\sum_{k=1}^\infty \mu_kXe_k d\beta_k,\quad t\ge 0,$
where $\{\beta_k\}$ is a sequence of  independent real Brownian motions on a filtered probability space $(\Omega,\mathcal F,\{\mathcal F_t\},\P)$ and $\{e_k\}$ is an orthormal  basis in $L^2(\mathcal O)$ which for convenience will be taken as the eigenfunction system for the Laplace operator with Dirichlet boundary conditions, i.e., $-\Delta e_k=\lambda_ke_k\;\mbox{\rm in}\; \mathcal O,\; e_k=0\;\mbox{\rm on}\; \partial\mathcal O,$
where $\mathcal O$ is an open and bounded subset  of $\R^d$, with smooth boundary $\partial\mathcal O$.
We shall assume that $\sum_{k=1}^\infty \mu_k^2\lambda^2_k<\infty.$
Equation (\ref{e1.1}) for $0<\alpha<1$ is relevant in the mathematical modelling of the dynamics of an ideal gas in a porous medium and, in particular, in a plasma fast diffusion model (for $\alpha=1/2$) (see e.g.  \cite{6}). The existence and uniqueness of a strong solution in the sense to be defined below was studied in  \cite{4},\cite{3},\cite{5},\cite{10} for more general nonlinear stochastic equations of the form (\ref{e1.1}).
In \cite{5} (see also \cite{4}) it was also proven  that for $\alpha=0$ and $d=1$ the solution $X=X(t,x)$ to (\ref{e1.1}) has the finite extinction property:
$
\P(\tau\le n)\ge 1-\frac{|x|_{-1}}{\rho\gamma}\;\left(\int_0^n e^{-C_N s}ds   \right) ^{-1}
$
for $|x|_{-1}<C_N^{-1}\rho\gamma$ where $\tau=\inf\{t\ge 0:|X(t,x)|_{-1}=0\}=\sup\{t\ge 0:|X(t,x)|_{-1}>0\}$ and $C_N, \gamma$ are constants related to the Wiener process $W$ and respectively to the domain $\mathcal O\subset \R^1$.

The following notations will be used in the sequel.
$H=L^2(\mathcal O),\;p\ge 1,$ with the norm denoted by $|\cdot|_2$ and  scalar product $\langle \cdot, \cdot   \rangle$.
$H^{-1}(\mathcal O)$ is the  dual of the  Sobolev space $H_0^{1}(\mathcal O)$ and is endowed with the scalar product $\langle u,v   \rangle_{-1}=\langle u,(-\Delta)^{-1}v\rangle$, where $\Delta$ is the Laplace operator with    domain $H^{2}(\mathcal O)\cap H_0^{1}(\mathcal O)$.
 All processes $X=X(t)$ arising here are adapted with  respect to the filtration $\{\mathcal F_t\}$.
For a Banach space $E$,  $L^p_W(0,T;E)$ denotes the  space of all   adapted processes in $L^p(0,T;E)$.
We shall use standard notation for Sobolev spaces  and spaces of integrable functions on $\mathcal O$.

\section{The main result}
\begin{definition}
\label{d2.1}
Let $x\in H$. An $H$-valued   continuous $(\mathcal F_t)$-adapted process $X=X(t,x)$ is called a solution to (\ref{e1.1}) on $[0,T]$ if 
$
X\in L^{p} (\Omega\times(0,T)\times \mathcal O)\cap L^2(0,T;L^2(\Omega,H)),\; p\ge 2,$  
  such that $\P$-a.s. $\forall\;j\in \N,\;t\in [0,T]$,
\begin{equation}
\label{e2.1}
\begin{array}{lll}
\langle X(t,x),e_j\rangle &=&\langle x,e_j\rangle +\int_0^t\int_\mathcal O(\rho |X(s,x)(\xi)|^{\alpha}\;\mbox{\rm sign}\;X(s,x)(\xi)+\tilde \Psi(X(s,x)(\xi)))\Delta e_j(\xi) d\xi ds\\
&&+\sum_{k=1}^\infty\mu_k\int_0^t\langle X(s,x)e_k,e_j\rangle d\beta_k(s),
\end{array}
\end{equation}
\end{definition}
For $x\in L^p(\mathcal O)$, $p\ge 4$ and $d=1,2,3$ there is a unique solution
$
X\in L_W^{\infty}(0,T;L^p(\Omega,H))
$
to (\ref{e1.1}) in the sense of Definition \ref{d2.1}. Moreover, if $x\ge 0$ a.e. in $\mathcal O$
then $X\ge 0$ a.e. in $\Omega\times[0,T]\times\mathcal O)$.

By the proof of \cite[Theorem 2.2]{5} and \cite[Proposition 3.4]{5} we also know  that for $\lambda\to 0$,
\begin{equation}
\label{e2.3}
\left\{\begin{array}{l}
X_\lambda\to X\quad\mbox{\rm strongly both in}\;L^2(0,T;L^2(\Omega,L^2(\mathcal  O)))\;
\mbox{\rm and in $L^2(\Omega;C([0,T];H))$},\\ \mbox{\rm  weakly in}\;L^p(\Omega\times(0,T)\times\mathcal  O),
\mbox{\rm and weak$^*$   in}\;L^\infty(0,T;L^p(\Omega;L^p(\mathcal  O))),
\end{array}\right. 
\end{equation}
 where $X_\lambda ,\;\lambda>0, $ is the solution to   approximating equation
\begin{equation}
\label{e2.4}
\left\{\begin{array}{l}
dX_\lambda(t)-\Delta(\Psi_\lambda(X_\lambda(t))+
\lambda X_\lambda(t)+\tilde\Psi(X_\lambda(t)))
dt= \sigma(X_\lambda(t)) dW(t), 
\\
\Psi_\lambda(X_\lambda)+\lambda X_\lambda+\tilde\Psi(X_\lambda)=0\quad\mbox{\rm on}\;\partial\mathcal  O,\quad X_\lambda(0,x)=x,\\
\Psi_\lambda(x)=\frac1\lambda\;(x-(1+\lambda\Psi_0)^{-1}(x))= \Psi_0((1+\lambda\Psi_0)^{-1}(x)),\quad \Psi_0(x)=\rho|x|^\alpha\;\mbox{\rm sign}\;x.
\end{array}\right.
\end{equation}
 Everywhere in the sequel $X=X(t,x)$ is the solution to (\ref{e1.1}) in the sense of Definition \ref{d2.1} where $x\in L^4(\mathcal  O)$.
Below  $\gamma$ shall denote the minimal constant arising in the Sobolev embedding $L^{\alpha+1}(\mathcal  O)\subset H^{-1}(\mathcal  O)$ (see  (\ref{e3.5})  below) and   $C^*=\sum_{k=1}^\infty\mu^2_k\;|e_k|^2_{H^{1}_0(\mathcal  O)}=\sum_{k=1}^\infty\mu^2_k\;\lambda_k^2.$ 
Theorem \ref{t2.1} is the main result of  the  paper.
\begin{theorem}
\label{t2.1}
Assume that $d=1,2,3$ and   that $
0<\alpha<1\;\mbox{\rm if $d=1,2$}, \;
\frac15\le\alpha<1\;\mbox{\rm if $d=3$}.$ Let $\tau:=\inf\{t\ge 0:\;|X(t,x)|_{-1}=0\}$.
Then we have $|X(t,x)|_{-1}=0,\;\mbox{\rm for $t\ge \tau$},\;\P\mbox{\rm -a.s.}.$
Furthermore $\P(\tau\le t)\ge 1-\frac{|x|^{1-\alpha}_{-1}}{(1-\alpha)\rho\gamma^{1+\alpha}}\;
 \left( \int_0^te^{-(1-\alpha) C^*s}ds\right)^{-1}.$
In particular, if $|x|^{1-\alpha}_{-1}< \frac{\rho\gamma^{1+\alpha}}{C^*}$, then $\P(\tau<\infty)>0,$ and if $C^*=0$, then $\tau\le \frac{|x|^{1-\alpha}_{-1}}{(1-\alpha)\rho\gamma^{1+\alpha}}$.
\end{theorem}
\begin{remark}
\em 
This result extends to  $\mathcal O\subset \R^d$ with $d\ge 4$, if $\alpha\in [\frac{d-2}{d+2},1)$. However, we have to strengthen the assumption on  $\mu_k$, $k\in \N$, see \cite[Section 4]{4} and in particular \cite[Remark 2.9(iii)]{8} for a detailed discussion.
 \end{remark}

\section{\bf Proof of Theorem \ref{t2.1}}
We shall proceed as in  the  proof of \cite[Theorem 4.2]{5}. Consider the solution $X_\lambda\in L^2_W(0,T;L^2(\Omega;H^1_0(\mathcal O)))$ to equation (\ref{e2.4}). Then by applying the classical It\^o formula  to the real valued semi-martingale $|X_\lambda(t)|^2_{-1}, t\in [0,T],$
and to the function
$
\varphi_\varepsilon(r)=(r+\varepsilon^2)^{(1-\alpha)/2},\quad r\in \R,
$
we find that
\begin{equation}
\label{e3.1}
\begin{array}{l}
d\varphi_\varepsilon(|X_\lambda(t)|^2_{-1})+(1-\alpha)(|X_\lambda(t)|^2_{-1}+\varepsilon^2)^{-(1+\alpha)/2}\langle X_\lambda(t),\Psi_\lambda(X_\lambda(t)) +\lambda X_\lambda(t)+  \tilde\Psi_\lambda(X_\lambda(t)) \rangle dt \\
=\frac12\;\sum_{k=1}^\infty\mu_k^2(1-\alpha)
\frac{|X_\lambda(t)e_k|^2_{-1}(|X_\lambda(t)|^2_{-1}+\varepsilon^2)-(1-\alpha)^2|\langle X_\lambda(t)e_k,X_\lambda(t)   \rangle_{-1}|^2)}{(|X_\lambda(t)|_{-1}^2+\varepsilon^2)^{(3+\alpha)/2}}\;dt\\
+\langle\sigma(X_\lambda(t))dW(t),\varphi'_\varepsilon(|X_\lambda(t)|^2_{-1})X_\lambda(t)   \rangle_{-1}\\
\le \frac12\;\sum_{k=1}^\infty\mu_k^2\;\frac{(1-\alpha)|X_\lambda(t)e_k|^2_{-1}}{(|X_\lambda(t)|_{-1}^2+\varepsilon^2)^{(1+\alpha)/2}}dt+\langle\sigma(X_\lambda(t))dW(t),\varphi'_\varepsilon(|X_\lambda(t)|^2_{-1})X_\lambda(t) \rangle_{-1}\\
\ds\le C^*\frac{(1-\alpha)|X_\lambda(t)e_k|^2_{-1}}{(|X_\lambda(t)|_{-1}^2+\varepsilon^2)^{(1+\alpha)/2}}dt+\langle\sigma(X_\lambda(t))dW(t),\varphi'_\varepsilon(|X_\lambda(t)|^2_{-1})X_\lambda(t) \rangle_{-1}
\end{array} 
\end{equation}
Then letting $\lambda\to 0$,  by (\ref{e2.3})
we get  that $\ds\liminf_{\lambda\to 0}\int_0^T\langle \Psi_{\lambda} (X_{\lambda}(t)),X_{\lambda}(t)\rangle dt$ $\ds\ge \rho\int_0^T|X(t)|^{1+\alpha}_{L^{1+\alpha}(\mathcal O)}dt$, $\P$-a.s. and hence  
\begin{equation}
\label{e3.2}
\begin{array}{l}
\ds\varphi_\varepsilon(|X(t)|^2_{-1})+ (1-\alpha)\rho\int_r^t\frac{|X(s)|^{\alpha+1}_{L^{\alpha+1}(\mathcal O)}}{(|X(s)|_{-1}^2+\varepsilon^2)^{(1+\alpha)/2}}ds\le \varphi_\varepsilon(|X(r)|^2_{-1})
\\
\ds+C^*\int_r^t\frac{(1-\alpha)|X(s)|^{2}_{-1}}{(|X(s)|_{-1}^2+\varepsilon^2)^{(1+\alpha)/2}}ds+2\int_r^{t}\langle\sigma(X(s))dW(s),\varphi'_\varepsilon(|X(s)|^2_{-1})X(s) \rangle_{-1},
\quad\P\mbox{\rm -a.s.},\;r<t.
\end{array} 
\end{equation}
Next by the Sobolev embedding theorem we have
\begin{equation}
\label{e3.5}
\begin{array}{l}
|u|_{-1}\le \gamma |u|_{L^{\alpha+1}(\mathcal O)},\;\forall\;u\in L^{\alpha+1}(\mathcal O),\;\mbox{\rm if $d>2$ and $\alpha\ge \frac{d-2}{d+2},$
and $\forall\;\alpha>0$, if d=1,2}.
\end{array}
\end{equation}
Then substituting (\ref{e3.5}) into (\ref{e3.2}) we get
\begin{equation}
\label{e3.2bis}
\begin{array}{l}
\ds\varphi_\varepsilon(|X(t)|^2_{-1})+ (1-\alpha)\rho\gamma^{1+\alpha}\int_r^t\frac{|X(s)|^{\alpha+1}_{-1}}{(|X(s)|_{-1}^2+\varepsilon^2)^{(1+\alpha)/2}}ds\le \varphi_\varepsilon(|X(r)|^2_{-1})
\\
\ds+C^*\int_r^t\frac{(1-\alpha)|X(s)|^{2}_{-1}}{(|X(s)|_{-1}^2+\varepsilon^2)^{(1+\alpha)/2}}ds+\int_r^{t}\langle\sigma(X(s))dW(s),\varphi'_\varepsilon(|X(s)|^2_{-1})X(s) \rangle_{-1},
\quad\P\mbox{\rm -a.s.},\;r<t.
\end{array} 
\end{equation}
 Now for $\epsilon\to  0$ we have
$$
\begin{array}{l}
\ds|X(t)|^{1-\alpha}_{-1}+ (1-\alpha)\rho\gamma^{1+\alpha}\int_r^t1_{\{|X(s)|_{-1}>0\}}ds\le |X(r)|^{1-\alpha}_{-1}+C^*(1-\alpha)\int_r^t|X(s)|^{1-\alpha}_{-1}ds\\
\ds+(1-\alpha)\int_r^{t}\langle\sigma(X(s))dW(s),|X(s)|^{-(\alpha+1)}_{-1}X(s) \rangle_{-1},
\;\P\mbox{\rm -a.s.},\;r<t.
\end{array} 
$$
Hence by It\^o's product rule
\begin{equation}
\label{e3.2ter}
\begin{array}{l}
\ds e^{-C^*(1-\alpha)t}|X(t)|^{1-\alpha}_{-1}+ (1-\alpha)\rho\gamma^{1+\alpha}\int_r^te^{-C^*(1-\alpha)s}1_{\{|X(s)|_{-1}>0\}}ds\le e^{-C^*(1-\alpha)r} |X(r)|^{1-\alpha}_{-1}\\
\ds+ (1-\alpha)\int_r^{t}e^{-C^*(1-\alpha)s}1_{\{|X(s)|_{-1}>0\}}\langle\sigma(X(s))dW(s),|X(s)|^{-(\alpha+1)}_{-1}X(s) \rangle_{-1},
\quad\P\mbox{\rm -a.s.},\;r<t.
\end{array} 
\end{equation}
From this it immediately follows that 
$e^{-C^*(1-\alpha)t}|X(t)|^{1-\alpha}_{-1},\;t\ge 0,$ is an $(\mathcal F_t)$-supermartingale, hence $|X(t)|_{-1}=0$ for all $t\ge \tau$. So, (\ref{e3.2ter}) with $r=0$ after taking expectation implies that
$
\int_0^{t}e^{-C^*(1-\alpha)s}\P(\tau>s)ds\le\frac{|x|^{1-\alpha}_{-1}}{(1-\alpha)\rho\gamma^{1+\alpha}},\quad t\ge 0.
$
This implies that
$
\P(\tau>t)\le\frac{|x|^{1-\alpha}_{-1}}{(1-\alpha)\rho\gamma^{1+\alpha}}\left(\int_0^{t}e^{-C^*(1-\alpha)s}ds\right)^{-1},\;t\ge 0,
$
and the assertion follows. $\Box$

\section*{Acknowledgement} 
 This work has been supported in part by
 the PIN-II ID-404 (2007-2010) project  of
Romanian Minister of Research,
the DFG -International Graduate School ``Stochastics and Real World Models'',the SFB-701 and the
BiBoS-Research Center.', 
  the research programme ``Equazioni di
Kolmogorov'' from the Italian
``Ministero della Ricerca Scientifica e Tecnologica''
and "FCT, POCTI-219, FEDER".

\bigskip

\end{document}